\documentclass{article}
\usepackage{amssymb,euscript,latexsym, amsmath, amsthm}
\usepackage[dvips]{graphicx}

\makeatletter
\@addtoreset{figure}{section}
\def\thefigure{\thesection.\@arabic\c@figure}
\def\fps@figure{h, t}
\@addtoreset{equation}{section}

\makeatother
\def\beq{\begin{equation}}
\def\eeq{\end{equation}}

\newtheorem{thm}{Theorem}[section]
\newtheorem{prop}[thm]{Proposition}
\newtheorem{lem}[thm]{Lemma}
\newtheorem{cor}[thm]{Corollary}
\newtheorem{dfn}[thm]{Definition}

\begin{document}

\title{Abstract mechanical connection and Abelian reconstruction for almost 
       K\"{a}hler manifolds.}

\author{Sergey Pekarsky,    Anthony D. Blaom, Jerrold E. Marsden \\
Control and Dynamical Systems 107-81 \\ California Institute of
Technology,\\ Pasadena, CA 91125}

\date{\small Preprint:  Summer 1998. 
Last revised: \today}
\maketitle

\begin{abstract}
When the phase space $P$ of a Hamiltonian $G$-system $(P, \omega, G, J, H)$
has an almost K\"{a}hler structure a preferred connection, called 
\emph{abstract mechanical connection}, 
can be defined by declaring horizontal spaces at each point 
to be metric orthogonal to the tangent to the group orbit.
Explicit formulas for the corresponding connection
one-form ${\mathcal{A}}$ are derived in terms of the momentum map, 
symplectic and complex structures.
Such connection can play the role of the \emph{reconstruction
connection} (see A. Blaom, 
Reconstruction phases via Poisson reduction.
{\em Diff. Geom. Appl.}, 12:231--252, 2000.),
thus significantly simplifying
computations of the corresponding dynamic and geometric phases
for an Abelian group $G$.
These ideas are illustrated using the example of the resonant three-wave 
interaction. Explicit formulas for the connection one-form and
the phases are given together with some new results on the
symmetry reduction of the Poisson structure.
\end{abstract}


\tableofcontents

\section{Introduction}

\paragraph{Definitions and preliminaries.}
Consider a finite dimensional symplectic manifold $(P, \omega)$.
Let a Lie group $G$ act on it canonically, {\it i.e.} by preserving the
symplectic form $\omega$, and assume that this action admits an
(equivariant) momentum map $J  : P \rightarrow U \subset 
{\mathfrak{g}}^\ast$, $U \equiv J(P)$.
Let a dynamical system be defined on $P$ by some Hamiltonian $H$.
We call $(P, \omega, G, J, H)$ a Hamiltonian $G$-system.
Assume also that $G$ acts on $P$ freely and properly so that the
Poisson reduction can be performed (in fact, these conditions can be 
slightly relaxed, see e.g. \cite{MaRa1994}).
For background on momentum maps, Poisson reduction, etc, the reader
is referred to Marsden and Ratiu \cite{MaRa1994}.

Recall that an almost K\"{a}hler manifold $(M, \omega, 
{\mathcal{J}}, {\mathfrak{s}})$ can be defined as a manifold $M$
with an almost complex structure ${\mathcal{J}}$ and an 
${\mathcal{J}}$-invariant (i.e. hermitian) metric ${\mathfrak{s}}$,
whose fundamental $2$-form $\omega$ defined by
\begin{equation}
  \label{akahler}
  \omega (x) ({\mathbf{v}},{\mathbf{w}}) = 
  {\mathfrak{s}} (x)({\mathcal{J}}(x){\mathbf{v}},{\mathbf{w}})
  \qquad \forall \  v, w \in T_x M
\end{equation}
is closed and hence a symplectic form on $M$. If in addition the
Nijenhuis torsion of ${\mathcal{J}}$ vanishes, then ${\mathcal{J}}$
is complex and $M$ becomes a  K\"{a}hler manifold \cite{KoNo1963}.
The automorphisms of an almost K\"{a}hler structure are
diffeomorphisms of $M$ which at the same time are symplectomorphisms,
almost complex maps and isometries w.r.t. 
$\omega, {\mathcal{J}}$ and ${\mathfrak{s}}$, resp.
It follows from the definition that any two of these conditions imply
the third one. For the background and more information see, e.g.
Kobayashi and Nomizu \cite{KoNo1963}.

\paragraph{Reconstruction of the dynamics.}
The space of group orbits $P/G$, which is obtained by taking the 
quotient map $\pi : P \rightarrow P/G$ and is a smooth manifold under
appropriate assumptions, inherits a Poisson structure from that of $P$.
The Hamiltonian $H$ drops to a reduced Hamiltonian $h$ on $P/G$, and
the corresponding Hamiltonian vector fields $X_H$ and $X_h$, as well as
their solutions $x_t$ and $y_y$, resp., are related by the projection 
$\pi : P \rightarrow P/G$. 

Assume that $y_t$ is periodic with period $T$, then for any initial 
condition $x_0 \in \pi^{-1}(y_0)$, the associated reconstruction phase
is the unique $g \in G$ such that $x_T = g \cdot x_0$. 
The methods presently used to compute reconstruction phases
are generally based on those established in \cite{MaMoRa90}.
The procedure can be sketched as follows.

If $J:P\rightarrow{\mathfrak g}^*$ is a momentum map, which 
we will suppose is $\operatorname{Ad}^*$-equivariant, then 
${J}(x_t)=\mu_0\equiv{J}(x_0)$.  
Under appropriate connectedness hypotheses  
the Marsden-Weinstein reduced space ${J}^{-1}(\mu_0)/G_{\mu_0}$ 
($G_{\mu_0}$ denoting the isotropy of the 
co-adjoint action at $\mu_0\in{\mathfrak g}^*$) can be identified with 
a symplectic leaf $P_{\mu_0}\subset P/G$ containing the reduced 
solution curve $ y_t$, and the  projection 
${J}^{-1}(\mu_0)\rightarrow P_{\mu_0}$ is a principal $G_{\mu_0}$-bundle.

The first step in calculating the reconstruction phase $ {g}$ above is 
to equip the bundle ${J}^{-1}(\mu_0)\rightarrow P_{\mu_0}$ 
with a principal connection $\alpha_{\mu_0}$, whose holonomy along the 
reduced curve $ y_t$ is called the associated \emph{geometric phase} and
denoted $g_{\operatorname{geom}}$. The phase $ {g}$ is then the product 
$g_{\operatorname{dyn}} g_{\operatorname{geom}}$ where 
$ g_{\operatorname{dyn}}$, called the \emph{dynamic 
phase}, is obtained by integrating a linear, non-autonomous, ordinary 
differential equation, called the \emph{reconstruction equation}.  The 
coefficients in this equation are defined in terms of 
$\alpha_{\mu_0}$, the unreduced Hamiltonian vector field $ X_H$, 
and an $\alpha_{\mu_0}$-horizontal lift of $y_t$ to ${\mathbf 
J}^{-1}(\mu_0)$.  Calculating the geometric phase usually requires one 
to compute the curvature of $\alpha_{\mu_0}$.

While any connection $\alpha_{\mu_0}$ can be used to compute $ {g}$ as 
above, a poor choice will lead to unwieldy computations.
For the so-called simple mechanical $G$-systems a natural choice exists 
(see below); for other systems the choice is often made on a case-by-case 
basis.

\paragraph{Overview of the results.}
As we already mentioned, the methods presently used to compute 
reconstruction phases
are based on those established in \cite{MaMoRa90}.
Though the general ideas in \cite{MaMoRa90} apply for arbitrary 
Hamiltonian $G$-systems, most of the advances in the computation
techniques have been done for mechanical systems on cotangent
bundles $T^* Q$ of some Riemannian manifolds $Q$ with the metric,
which determines the kinetic energy, playing the crucial role
for the definition of the \emph{mechanical connection}.
Unfortunately, these settings exclude such interesting and
important systems as $N$ point vortices on a plane or on a sphere,
$N$-wave interaction, etc. where the configuration space
is \emph{not} a cotangent bundle. 

Luckily, some of these systems have a natural almost K\"{a}hler
structure which we exploit in the construction of the
\emph{abstract mechanical connection}. It is defined by specifying
the horizontal space to be \emph{metric orthogonal} to the group
orbits (see Section $4$ for the details). The corresponding
connection one-form is then obtained in terms of the momentum map, 
symplectic and complex structures. These expressions enable
to further simplify the computation of the reconstruction phases
as described in \cite{Blaom99} for the case of Abelian groups.
The requirement of the group being Abelian is essential for the
geometric phase part (see Section $4.2$), but it is not used in the 
construction of the map $L$ involving the abstract inertia tensor
and in the expression for the dynamic phase. The Abelian
property of the group makes the relation between the principal 
connections on Poisson and symplectic bundles trivial (see Section $2.3$)
It also  significantly simplifies the picture of dual pairs, which
underlies our constructions, and enables us to construct a very
'useful' bundle $j : P/G \rightarrow {\mathfrak g}^*$. This is
briefly described in Section $2.4$.

In the work in progress a generalization to non-Abelian group action
is being considered, as well as further simplifications and links
which arise in the case of the phase space being 
a cotangent bundle with the almost K\"{a}hler structure coming
from a Riemannian metric on the configuration space.
In particular,  the relation 
between the abstract mechanical connection and the 
well-known mechanical connection are considered in \cite{Blaom99}, where
the reconstruction phases for the cotangent bundles were analyzed,
though not from the point of view
of almost K\"{a}hler manifolds and corresponding abstract
mechanical connections.

\section{Reconstruction connection and the associated phases}

In  this section we briefly overview the results on reconstruction
phases obtained in \cite{Blaom99}. We refer the reader to the
original paper for detailed and comprehensive treatment of the subject.
Here we are mainly interested in adopting these results to the case
of almost K\"{a}hler systems, and thus we shall avoid giving much
details to keep the presentation clear and avoid repetition.

In \cite{Blaom99} a general formula is derived which expresses a 
reconstruction phase in terms of the associated reduced solution,
viewed as a curve in the Poisson-reduced phase space $P/G$, and certain
derivatives \emph{transverse} to the symplectic leaf in $P/G$ containing 
the curve. Specifically, the dynamic part of the phase depends on 
transverse derivative in the Poisson-reduced Hamiltonian, while
the geometric part is determined by transverse derivatives in the
leaf symplectic structures.

\paragraph{Highlights and basic assumptions.}
It is shown in \cite{Blaom99} that the principal connection on the bundle
$J^{-1} (\mu_0) \rightarrow P_{\mu_0}$, which plays a crucial
role in the computation of the phases, is most naturally
viewed as the restriction to $J^{-1} (\mu_0)$ of a certain kind of
distribution $A$ on $P$, which is called a \emph{reconstruction
connection}. To define the transverse derivatives above one then 
specifies a connection $D$ on the symplectic stratification of $P/G$ 
(a distribution on $P/G$ furnishing a complement for the characteristic
distribution). This connection $D$ can be obtained by 'Poisson-reducing'
the connection $A$.


Explicitly, assuming that as a cycle the reduced  curve $y_t$ above 
is a boundary $\partial\Sigma$ ($\Sigma\subset P_{\mu_0}$ compact and 
oriented), the corresponding reconstruction  phase is 
$g=g_{\operatorname{dyn}}
g_{\operatorname{geom}}$, where
$$
g_{\operatorname{dyn}} =\exp
\int_0^T D_{\mu_0}h(y_t)\,dt
$$
$$
g_{\operatorname{geom}} =\exp
\int_\Sigma  D_{\mu_0}\omega_D
$$

In these formulas $D_{\mu_0}$ denotes a certain `exterior covariant 
derivative' depending on $ D$ and $\mu_0$ that maps ${\mathbb R}$-valued 
$p$-forms on $ P/G$ to ${\mathfrak g}_{\mu_0}$-valued 
$p$-forms on $ P_{\mu_0}$ ($p=0, 1, 2,\ldots$).  For example, 
$D_{\mu_0}h(y_t)$ is an element of ${\mathfrak g}_{\mu_0}$ that happens 
to measure the derivative of $ h$ in directions lying in $ D(y_t)$ 
(i.e., in certain directions transverse to the symplectic leaves).  
$\omega_D$ denotes the two-form on $ P/G$ whose 
restriction to a given leaf gives that leaf's symplectic structure, 
and whose contraction with vectors in $ D$ vanishes .

These formulas apply assuming that $G_{\mu_0}$ is 
Abelian, that $ P_{\mu_0}$ is a non-degenerate symplectic leaf, and 
that $ D$ is a smooth distribution in some neighborhood of $ P/G$.  
These conditions are in addition to the following assumptions
which are understood to be in place throughout the paper:
\begin{itemize}
\item all manifolds are  smooth, i.e. $C^\infty$
\item the group $G$ acts freely and properly, so that the natural
projection $\pi : P \rightarrow P/G$ is a submersion
\item the group action is \emph{Hamiltonian}, i.e. it admits
a momentum map $J : P \rightarrow {\mathfrak g}^*$ which is
$\operatorname{Ad}^*$-equivariant: $J (g \cdot x) = g \cdot (x)$.
\end{itemize}

The formulas above make sense for \emph{any} connection $D$ on the 
symplectic stratification of $ P/G$; whence, the total reconstruction 
phase $g=g_{\operatorname{dyn}} g_{\operatorname{geom}}$ 
(which is independent of the choice of 
$\alpha_{\mu_0}$, and hence $ A$) can be computed  using {\em  any} 
connection $ D$ on the symplectic stratification of $ P/G$.

The following two subsections give a short review of the results in
\cite{Blaom99} that are relevant for our applications.

\subsection{Main constructions.}

\begin{dfn}
Call a distribution $A$ on $P$ a \emph{reconstruction connection} if
\begin{itemize}
\item[a)] $A$ is $G$-invariant
\item[b)] $\operatorname{Ker} T_x J = T_x (G_\mu \cdot x) \oplus A(x)
\quad (x \in P, m \equiv J(x))$
\end{itemize}
Here $G_\mu$ denotes the point stabilizer of the coadjoint action at 
$\mu \in {\mathfrak g}^*$, $T_x J$ is the tangent map, and $\oplus$
denotes the direct sum.
\end{dfn}

\paragraph{Connections on the symplectic stratification of $P/G$.}
Let $E$ denote the characteristic distribution on $P/G$ (i.e. the
distribution tangent to the symplectic leaves). We call a distribution $D$
on $P/G$ a \emph{connection on the symplectic stratification of $P/G$}
if it furnishes a complement for $E$:
\begin{equation}
\label{distr_D}
T(P/G) = E \oplus D.
\end{equation}

Now let $A$ be a $G$-invariant distribution on $P$. Since 
$G$ acts by symplectic diffeomorphisms, the distribution $A^\omega$
is also $G$-invariant. It consequently drops to a distribution
$\hat{A} \equiv \pi_* (A^\omega)$ on $P/G$; here $\pi_*$ denotes
push-forward. Conversely, if $D$ is an arbitrary distribution on $P/G$,
then $\hat{D} \equiv (\pi^* D)^\omega$ is a $G$-invariant distribution
on $P$; here $\pi^*$ denotes pull-back. Evidently, one has 
$$
\hat{\hat{D}} = D.
$$

We quote the following theorem form \cite{Blaom99} without proof.
\begin{thm}[Blaom 1999]
\label{thm_A}
If $A$ is a general reconstruction connection, then $\hat{A}$ is
a connection on the symplectic stratification of $P/G$.
Moreover, the map $A \mapsto \hat{A}$ is a bijection from the set
of reconstruction connections to the set of connections on the
symplectic stratification of $P/G$. This bijection has an inverse
$D \mapsto \hat{D}$.
\end{thm}

If $A$ is a reconstruction connection, one thinks of $\hat{A}$ as its 
Poisson-reduced counterpart. A reconstruction connection $A$ can be
reconstructed from its reduced counterpart $D \equiv \hat{A}$ according
to $A = \hat{D}$.

Two other lemmas from \cite{Blaom99} are relevant to our presentation
and will be used in Section $3$.
\begin{lem}
\label{lem_B1}
Let $\pi : P \rightarrow Q$ be a Poisson submersion and let $E$
denote the characteristic distribution on $Q$. If $P$ is symplectic,
and $\omega$ denotes the symplectic form on $P$, then
$$
\pi^* E = \operatorname{Ker} T \pi + (\operatorname{Ker} T \pi)^\omega.
$$
\end{lem}

\begin{lem}
\label{lem_B2}
Let $x \in P$ be arbitrary and define $\mu \equiv J(x)$. Then
$$
T_x (G_\mu \cdot x) = ((\pi^*E)(x))^\omega.
$$
\end{lem}

\paragraph{Transverse derivatives in $P/G$.}
Under appropriate connectedness hypotheses each 
reduced space $J^{-1}(\mu)/G_\mu$ may be identified with a symplectic
leaf $P_\mu \subset P/G$. A connection $D$ on the symplectic 
stratification of $P/G$ allows one to define derivatives of
functions on $P/G$ transverse to $P_\mu$. At a point in $P_\mu$
such a derivative can be identified in a natural way with an
element of the isotropy algebra ${\mathfrak g}_\mu$, provided
the isotropy group $G_\mu$ is Abelian. More generally, for
such $\mu$ the connection $D$ defines an 'exterior 
covariant derivatives' mapping ${\mathbb R}$-valued $p$-forms on
$P/G$ to ${\mathfrak g}_\mu$-valued $p$-forms on the leaf $P_\mu$.

Let $D$ be a fixed connection on the symplectic stratification of $P/G$.
Fix $\mu \in U \equiv J(P)$ and assume $G_\mu$ is Abelian. Then:
\begin{prop}
\label{lem_L}
For each $y \in P_\mu$ there is a natural isomorphism 
$D(y) \leftrightarrow {\mathfrak g}_\mu^*$ well defined by
\begin{equation}
\label{L_map}
v \mapsto p_\mu ( \operatorname{forg} (T J \cdot w)),
\end{equation}
where $w$ denotes any element of $T_x P$ with $T \pi \cdot w = v$, and
$x \in J^{-1}(\mu) \cap \pi^{-1}(y) = \pi_\mu^{-1} (y)$ is arbitrary.
The inverse of this map (which depends on $D, \mu$ and $y$) is denoted
by $L(D, \mu, y) : {\mathfrak g}_\mu^* \rightarrow D(y)$.
\end{prop}

The map forg : $T U \rightarrow {\mathfrak g}^*$ denotes the map that
'forgets base point' and 
$p_\mu :{\mathfrak g}^* \rightarrow {\mathfrak g}_\mu^*$ denotes the natural
projection.

\begin{dfn}
Suppose $f$ is a function on $P/G$ defined in some neighborhood of $y$.
Then the \emph{$(D, \mu)$-exterior covariant derivative} of $f$ at $y$,
denoted $D_\mu f (y) \in {\mathfrak g}_\mu$, is defined through
\begin{equation}
\label{ext_der1}
\langle \nu , D_\mu f (y) \rangle = \langle d f, L(D, \mu, y) (\nu) \rangle
\qquad \forall \ \nu \in {\mathfrak g}_\mu^*.
\end{equation}
\end{dfn}

\begin{dfn}
Let $\sigma$ be a differential $p$-form on $P/G$ defined in a neighborhood
of $P_\mu$, and assume that $G_\mu$ is Abelian. Then the 
\emph{$(D, \mu)$-exterior covariant derivative} $D_\mu \sigma$ of $\sigma$
is the ${\mathfrak g}_\mu$-valued $p$-form on $P_\mu$ defined through
\begin{equation}
\label{ext_der2}
\langle \nu , D_\mu \sigma (v_1, \dots , y_p) \rangle = 
d \sigma ( L(D, \mu, y) (\nu), v_1, \dots , v_p )
\end{equation}
where $\nu \in {\mathfrak g}_\mu^*, v_1, \dots , v_p \in T_y P_\mu$
and $y \in P_\mu$.
\end{dfn}

\paragraph{Smoothness conditions.}
Let $A$ be a reconstruction connection and let $D$ be a connection on the
symplectic stratification of $P/G$. Then we say that $A$ is 
\emph{$\mu$-smooth} ($\mu \in U$) if the set
$$
\{ A(x) \ | \ x \in J^{-1} (\mu) \}
$$
is a smooth sub-bundle of the tangent bundle $T (J^{-1} (\mu))$.
We call $D$ \emph{$\mu$-smooth} if the set
$$
\{ D(y) \ | \ y \in P_\mu \}
$$
is a smooth sub-bundle of 
$T_{P_\mu} (P/G) \equiv \{ T_y (P/G) | y \in P_\mu \}$.

Then, the following smoothness results hold \cite{Blaom99}
\begin{itemize}
\item $D$ is $\mu$-smooth if and only if $A$ is $\mu$-smooth.
\item If $D$ is $\mu$-smooth, then $L(D,\mu,y)$ in (\ref{L_map}) 
      depends smoothly on $y \in P_\mu$.
\item If  $D$ is $\mu$-smooth
      then $D_\mu f : P_\mu \rightarrow {\mathfrak g}_\mu$ is smooth.
\item Similarly, for a $p$-form $\sigma$, $\mu$-smoothness of $D$
      ensures smoothness of $D_\mu \sigma$.
\end{itemize}

\subsection{Reconstruction phases.}

Let $H$ be a $G$-invariant Hamiltonian on $P$ and $h : P/G \rightarrow
{\mathbb R}$ its Poisson-reduced counterpart. 
With the assumptions stated in the previous section satisfied, 
consider an integral curve $x_t \in P$ of $X_H$. The curve remains
in the submanifold $J^{-1}(\mu_0) (\mu_0 \equiv J (x_0))$ for all
time $t$ for which it is defined.The Marsden-Weinstein reduction
bundle 
$$
\pi_{\mu_0} : J^{-1}(\mu_0) \rightarrow P_{\mu_0}
$$
is a principal $G_{\mu_0}$-bundle. Let
$$
\alpha_{\mu_0} : T (J^{-1}(\mu_0)) \rightarrow {\mathfrak g}_{\mu_0}
$$
denote the connection one-form on this bundle whose associated
horizontal space at each $x \in J^{-1}(\mu_0)$ is 
$\operatorname{hor}_{\mu_0} \equiv A(x)$.
To ensure that $\alpha_{\mu_0}$ and $\operatorname{hor}_{\mu_0}$ are
smooth, we require that $A$ be $\mu_0$-smooth.

Let $y_t \in P_{\mu_0}$ denote the integral curve of the reduced
Hamiltonian vector field $X_h$ on $P/G$ that has $y_0 = \pi (x_0) \in 
P_{\mu_0}$ as its initial point. Then as $X_H$ and $X_h$ are $\pi$-related,
we have $y_t = \pi (x_t)$ for all $t$.

Let $d_t \in J^{-1}(\mu_0)$ denote the 
$\operatorname{hor}_{\mu_0}$-horizontal lift of $y_t$ having $x_0$ 
as its initial point $d_0$. Supposing that $y_t$ is periodic with
period $T$, we have 
$$
d_T = g_{\operatorname{geom}} \cdot x_0
\qquad
x_T = g_{\operatorname{dyn}} \cdot d_T,
$$
for some uniquely defined $g_{\operatorname{geom}}, g_{\operatorname{dyn}}
\in G_{\mu_0}$ called \emph{geometric} and \emph{dynamic} phases
associated with the reduced solution $Y_t$. 
The product $g_{\operatorname{total}} = g_{\operatorname{geom}}
 g_{\operatorname{dyn}}$ is called the \emph{total} phase. It does not
depend on $A = \hat{D}$, but depends only on $y_0$, the flow of $X_H$,
and the period $T$.

\paragraph{Dynamic phases.}
It is well known (\cite{MaMoRa90}) that the dynamic phase is given
by the solution of the following initial value problem,
known as the \emph{reconstruction equation}:
$$
\dot{g}_t = g_t \xi_t, \quad \operatorname{where} \quad
\xi_t \equiv \alpha_{\mu_0} (X_H (d_t)), \ 
\operatorname{and} \ g_0 = \operatorname{Id}.
$$
Here $g_t \xi_t$ denotes the tangent action of $g_t$.

Corollary $3.6$ of \cite{Blaom99} states that assuming $G_{\mu_0}$
is Abelian, the dynamic phase is given by
\begin{equation}
\label{dyn_ph}
g_{\operatorname{dyn}} = \exp \int_0^T D_{\mu_0} h (y_t) dt.
\end{equation}

\paragraph{Geometric phases.}
Recall that the geometric phase $g_{\operatorname{geom}}$ associated
with a solution $x_t$ is the holonomy of a principal connection 
$\alpha_{\mu_0}$ on $J^{-1}(\mu_0) \rightarrow P_{\mu_0}$ along
the corresponding reduced solution curve $y_t = \pi (x_t) \in P_{\mu_0}$.
Assuming $G_{\mu_0}$ is Abelian, the holonomy of appropriate curves
is determined by the curvature of $\alpha_{\mu_0}$.
It is well known (\cite{MaMoRa90}) that if the cycle $y_t$ is
in fact a boundary $\partial \Sigma$ ($\Sigma \subset P_{\mu_0}$
compact and oriented), then 
\begin{equation}
\label{geom_ph_curv}
g_{\operatorname{geom}} = \exp \left( - \int_\Sigma \Omega_{\mu_0} \right),
\end{equation}
where $\Omega_{\mu_0}$ is the curvature of $\alpha_{\mu_0}$, viewed
as a ${\mathfrak g}_{\mu_0}$-valued two-form on the reduced space
$P_{\mu_0}$.

Theorem C of \cite{Blaom99} shows that \emph{all} curvature
information on $\alpha_{\mu_0}$ is encoded in: (i) the connection $D$
on the symplectic stratification of $P/G$ corresponding to the
reconstruction connection $A$, together with (ii) the Poisson
structure on $P/G$.

The connection $D$ allows to 'assemble' the reduced symplectic
structures $\omega_\Lambda$ ($\Lambda \subset P/G$ a symplectic leaf)
into a single two-form $\omega_D$ on $P/G$ by decreeing that
\begin{equation}
\label{omegaD}
\omega_D (u, v) \equiv \omega_\Lambda (p_D u, p_D v), \quad
u,v \in T (P/G),
\end{equation}
where $\Lambda$ denotes the leaf to which the common base point
of $u$ and $v$ belongs, and where $p_D : T (P/G) \rightarrow E$
denotes the projection along $D$ onto the characteristic distribution
$E$.

We remark that in general $\omega_D$ need not be smooth, but, if
$P_{\mu_0}$ is a non-degenerate symplectic leaf, then $\omega_D$ 
is smooth wherever $D$ is of constant rank and smooth.
Then, Corollary $4.5$ of \cite{Blaom99} states that assuming $G_{\mu_0}$
is Abelian and $\omega_D$ is smooth in a neighborhood of $P_{\mu_0}$,
the geometric phase is given by
\begin{equation}
\label{geom_ph}
g_{\operatorname{geom}} = \exp \int_\Sigma D_{\mu_0} \omega_D.
\end{equation}

\section{Connections on various bundles for Abelian groups.}

In  this section the relation between connections on
Poisson and symplectic bundles is analyzed. This establishes the validity 
of the application of results in \cite{Blaom99} to our settings
in the case of Abelian groups $G$,
so that the metric orthogonal spaces to the group orbit in the
whole tangent $T_x P$ as well as within the kernel 
$\operatorname{Ker} T J(x) \subset T_x P$ both constitute valid
horizontal spaces for Poisson and symplectic bundles, respectively.

In the final subsection, the formalism of dual pairs is introduced
into the picture. The Symplectic Leaf Correspondence Theorem 
brings insight into the structure of various bundles and relates
the corresponding connections. 
For the Abelian case, it gives a new interpretation of the connection
on symplectic stratification $D$ as a connection on the bundle
$j : P/G \rightarrow U \subset {\mathfrak g}^*$ of symplectic leaves
over the dual of the Lie algebra (see Section $2.4$).

\subsection{Connections on Poisson and symplectic bundles.}

Consider the relation between a connection
on the Poisson reduction bundle $P \rightarrow P/G$ and 
connections on each of the symplectic Marsden-Weinstein reduction
bundles $J^{-1} (\mu) \rightarrow J^{-1} (\mu)/G_\mu$
for different $\mu \in {\mathfrak{g}}^\ast$.
This relation can be easily established in the case of an Abelian
group $G$ when $G_\mu \equiv G$ and both bundles have similar fibers.

Recall that a connection on the bundle $P \rightarrow P/G$ is a 
Lie algebra valued one-form ${\mathcal A}$ on $P$ that is $G$-equivariant
$g \cdot {\mathcal A} = \operatorname{Ad}_g \cdot {\mathcal A}$ and
satisfies ${\mathcal A} (\xi_P) = \xi \ \forall \ \xi \in {\mathfrak g}$.
The corresponding horizontal space is defined by 
$\operatorname{hor} = \operatorname{Ker} {\mathcal A}$.
The following theorem then holds.

\begin{thm}
\label{thm_bundle}
For the case of an Abelian group $G$, a connection ${\mathcal A}$
on the Poisson bundle  induces connections $\alpha_\mu$ on symplectic
Marsden-Weinstein bundles for regular momentum values $\mu$.
In particular, it defines a \emph{reconstruction connection} $A$ on $P$. 
Moreover, the connections on the symplectic stratification
of $P/G$ corresponding to $A$
and to ${\mathcal A}$ coincide, {\it i.e.}
\begin{equation}
\label{D=D'}
\hat{A} = D = \hat{{\mathcal A}}
\end{equation}
\end{thm}

\begin{proof}
Choose a regular value $\mu \in {\mathfrak{g}}^\ast$
such that the symplectic reduction at $\mu$ is defined.
Define induced horizontal and vertical spaces at $x \in J^{-1}(\mu)$
by the intersections with $\operatorname{Ker} T J$:
$$
\operatorname{hor}_\mu=\operatorname{hor} \cap \operatorname{Ker} T J,
\qquad
\operatorname{ver}_\mu=\operatorname{ver} \cap \operatorname{Ker} T J.
$$
By definition, 
$\operatorname{hor}_\mu \cap \operatorname{ver}_\mu = 0$.
As $G$ is Abelian, $G_\mu = G$ and 
$\operatorname{Ker} T \pi \subset (\operatorname{Ker} T \pi)^\omega
=\operatorname{Ker} T J$, so that 
$\operatorname{ver}_\mu = \operatorname{Ker} T \pi$.
Using the following set-theoretical identity
$(A + B) \cap C = A + B \cap C$ if $A \subset C$, we obtain 
\begin{multline}
\operatorname{Ker} T J  = 
(\operatorname{Ker} T \pi + \operatorname{hor})
\cap \operatorname{Ker} T J \\
 = \operatorname{Ker} T \pi + \operatorname{hor} \cap \operatorname{Ker} T J
 = \operatorname{ver}_\mu + \operatorname{hor}_\mu.
\end{multline}
Hence, $\operatorname{Ker} T J = 
\operatorname{ver}_\mu \oplus \operatorname{hor}_\mu$.
The corresponding connection one-form $\alpha_\mu$ is defined
by the horizontal space via 
$\operatorname{Ker} \alpha_\mu = \operatorname{hor}_\mu$.
The collection of these $\alpha_\mu$ then define a reconstruction
connection $A$ as defined in Section $2.1$.
It is $G$-invariant because ${\mathcal A}$ is $G$-invariant
for Abelian groups.

Finally, for the connections on the symplectic stratification
of $P/G$ determined by connections on Poisson  and symplectic bundles, 
{\it i.e.} by ${\mathcal A}$ and $A$, resp. we have at $y = \pi(x)$:
$$
D'(y) \equiv \hat{{\mathcal A}}(x) = T \pi (\operatorname{hor}^\omega(x)) 
$$
and
\begin{multline}
D(y) \equiv \hat{A}(x) = T \pi ((\operatorname{hor_\mu})^\omega)
= T \pi ((\operatorname{hor} \cap \operatorname{Ker} TJ)^\omega) \\
= T \pi ((\operatorname{hor})^\omega + (\operatorname{Ker} TJ)^\omega) 
= T \pi ((\operatorname{hor})^\omega + \operatorname{Ker} T \pi)
= T \pi (\operatorname{hor}^\omega(x)),
\end{multline}
where $x \in J^{-1}(\mu)$ with $y=\pi (x)$ and
we have used that $(\operatorname{Ker} TJ)^\omega = 
\operatorname{Ker} T \pi$.
Comparing the last two expressions we conclude that $D = D'$.
\end{proof}

This result enables us to go back and forth between connections
on Poisson and symplectic bundles for Abelian groups; in particular,
it will let us apply results of \cite{Blaom99} for the reconstruction
phases and use the abstract mechanical connection (defined in the next
section) as a reconstruction connection.

\subsection{Connections on dual pairs}

Recall the notion of dual pairs introduced by Weinstein \cite{Weinstein83}.
Consider a symplectic manifold $(P, \omega)$, Poisson manifolds
$Q_1, Q_2$, and Poisson maps $\rho_i : P \rightarrow Q_i, \, i=1,2$.
If for almost all $x \in P$, $(\mbox{Ker} \ T \rho_1(x))^\omega =
\mbox{Ker} \ T \rho_2(x)$, the diagram 
$Q_1 \stackrel{\rho_1}{\longleftarrow} P
\stackrel{\rho_2}{\longrightarrow} Q_2$ is called a {\it dual pair}. 
The dual pair is called {\it full}, if $\rho_1, \rho_2$ are surjective
submersions. If $Q_1 \stackrel{\rho_1}{\longleftarrow} P
\stackrel{\rho_2}{\longrightarrow} Q_2$ is a full dual pair, then the spaces
of Casimir functions on $Q_1$ and $Q_2$ are in bijective correspondence,
i.e. ${\rm Cas}(Q_1) \circ \rho_1 = {\rm Cas}(Q_2) \circ \rho_2$ (Weinstein
\cite{Weinstein83}).

It was shown in  Adam and Ratiu \cite{AdRa88} 
that for a symplectic manifold
$(P, \omega)$ with a Hamiltonian action of a Lie group $G$ having an
equivariant momentum map $J  : P \rightarrow U \subset {\mathfrak{g}}^\ast$,
$U \equiv J(P)$, such that $\pi  : P \rightarrow P/G$ and $J$ are 
surjective submersions, $P/G \stackrel{\pi}\longleftarrow P 
\stackrel{J}\longrightarrow U$
is a full dual pair.
The Poisson reduced space $P/G$, being a base of a principle
$G$-bundle, is itself foliated by  symplectic leaves ${\Sigma}_y$
through points $y \in P/G$. 
Let us denote the space of symplectic leaves by ${\mathcal{S}}$.
With the proper connectedness assumptions these leaves are
precisely the symplectic
reduced spaces $P_\mu = J^{-1} (\mu)/G_\mu$ (note that
$G_\mu$ can be different for different values of $\mu$).

On the other hand, $P$ is foliated by the level sets of the momentum
map $J^{-1}(\mu)$, for different $\mu \in {\mathfrak{g}}^\ast$,
with the dual of the Lie algebra itself being a foliation
by coadjoint orbits  ${\mathcal{O}}_\mu$ through $\mu$.
It follows from the Symplectic Leaf Correspondence Theorem 
\cite{Weinstein83}
that (under the above assumptions) the base space
of this foliation is in one-to-one correspondence with ${\mathcal{S}}$,
the space of symplectic leaves of the Poisson reduced space $P/G$.
A natural one-to-one correspondence 
between the symplectic leaves in each leg of a dual pair has 
been described in Weinstein \cite{Weinstein83}, together with a sketch
of the proof.
Here we state the Symplectic Leaf Correspondence Theorem and refer
for a detailed and comprehensive proof to Blaom \cite{Blaom1998}.
\begin{thm}
Let $P$ be a symplectic manifold and 
$Q_1 \stackrel{\rho_1}\longleftarrow P \stackrel{\rho_2}\longrightarrow Q_2$
a full dual pair. Assume that each leg $\rho_j : P \rightarrow Q_j, \
j = 1,2$ satisfies the property that pre-images of connected sets are 
connected. Let ${\mathcal{F}}_j$ denote the set of symplectic leaves in
$Q_j$. Then, under the above assumptions, there exists a bijection 
${\mathcal{F}}_1 \rightarrow {\mathcal{F}}_2$ given by
$$
\Sigma_1 \mapsto \rho_2 (\rho_1^{-1}(\Sigma_1))
$$
having inverse
$$
\Sigma_2 \mapsto \rho_1 (\rho_2^{-1}(\Sigma_2)).
$$
\end{thm}

This theorem enables us to define a leaf-to-leaf bijection 
that maps symplectic leaves $\Sigma_y$ (which are diffeomorphic to
symplectic reduced spaces $P_\mu$) to coadjoint orbits 
${\mathcal{O}}_\mu$ in the dual of the Lie algebra, $\mu = J (x)$.
Yet another realization of the symplectic leaves $\Sigma_y$
is given by the Orbit Reduction Theorem \cite{MaRa1998} which
establishes one-to-one correspondence between orbit reduced spaces
$P_{{\mathcal{O}}_\mu} = J^{-1} ({\mathcal{O}}_\mu)/G$ 
and symplectic reduced spaces $P_\mu = J^{-1} (\mu)/G_\mu$.

In the case of an Abelian group $G$, the coadjoint orbits are
trivial, {\it i.e.} ${\mathcal{O}}_\mu = \{ \mu \}$, and 
$G_\mu = G$, so that $P_{{\mathcal{O}}_\mu} \equiv P_\mu$
and ${\mathcal{S}} \cong  {\mathfrak{g}}^\ast$.
It follows then from the Reduction Lemma (see, e.g. \cite{MaRa1998})
that $G$-orbits of any point $x \in P$ are \emph{isotropic}, {\it i.e.},
$T_x (G \cdot x) \subset (T_x (G \cdot x))^\omega$ or, equivalently,
$\operatorname{Ker} T_x \pi \subset (\operatorname{Ker} T_x \pi)^\omega$.
Moreover, the bijection ${\mathcal{F}}_1 \rightarrow {\mathcal{F}}_2$
becomes a well-define map of the manifolds,
$j : P/G \rightarrow U$ which can be obtained through factoring the
momentum map. Indeed, the equivariance of the momentum map $J : P
\rightarrow U \subset {\mathfrak g}^*$ amounts in the Abelian case to 
invariance. It therefore factors through $\pi : P \rightarrow P/G$,
delivering a map $j : P/G \rightarrow U$ making the following
diagram commute
\begin{figure}
\begin{center}
\includegraphics[width=4cm]{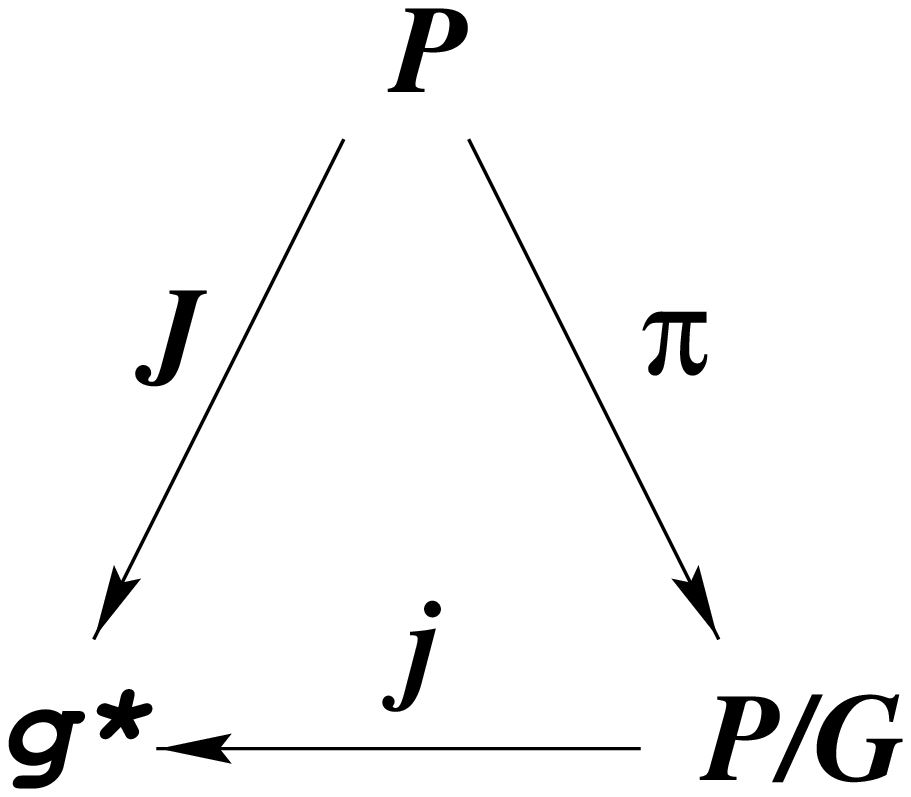}
\caption{The momentum map $J$  factors through
delivering a map $j : P/G \rightarrow U \subset {\mathfrak g}^*$.}
\label{diag_j}
\end{center}
\end{figure}

\vspace{2cm}

The map $j$ is a submersion since $J$ is a submersion (under our
hypothesis of a free action). Since the coadjoint orbits are points,
the symplectic leaves in $P/G$ are simply the fibers of $j$, that is,
$P_\mu = j^{-1} (\mu), \ \mu \in U$.

Thus, with this interpretation, the connection on the symplectic
stratification $D$ can be thought of as an (Ehresmann) connection on 
the bundle $j : P/G \rightarrow U \subset {\mathfrak g}^*$.
Theorem $2.1$ taken from \cite{Blaom99} as well as the results of the
previous subsection establish then a relation
between the connections $A$ and $D$ on the bundles 
$\pi : P \rightarrow P/G$ and $j : P/G \rightarrow U 
\subset {\mathfrak g}^*$, respectively.

Finally, the tangent map $T j$ delivers the isomorphism of the
Proposition \ref{lem_L}, where now $L$ does not depend 
explicitly on $\mu$ as $G$ is Abelian; 
that is, ${\mathfrak g}^*_\mu = {\mathfrak g}^*$,
$p_\mu \equiv \operatorname{Id}$ and the dependence on $\mu$ enters
only through $\mu = j (y)$.

\begin{lem}
\label{lem_L_Tj}
Let $L(D, y) : {\mathfrak g}^* \rightarrow D(y)$ be defined
by (\ref{L_map}), then its inverse is given by the tangent map $T j$
restricted to the distribution $D$:
\begin{equation}
\label{L_Tj}
L^{-1} = \operatorname{forg} \left( T j |_D \right) : D  
\rightarrow {\mathfrak g}^* ,
\end{equation}
where, the map forg : $T U \rightarrow {\mathfrak g}^*$ denotes the map that
'forgets base point'.
\end{lem}

\begin{proof}
The proof readily follows from the fact that the momentum map factors
through the quotient map, so that $T J = T j \circ T \pi$, and the
definition of the map $L$ for any $y \in P_\mu$ 
given by (\ref{L_map}), where $w$ is any vector 
in $T_x P$ that satisfies $T \pi \cdot w = v$, with $v \in D(y)$
and $x \in J^{-1}(\mu) \cap \pi^{-1} (x)$:
$$
v \mapsto L^{-1}(D, y) \cdot v \equiv p_\mu (\operatorname{forg}(TJ \cdot w))
= \operatorname{forg}(Tj \circ T \pi \cdot w) 
= \operatorname{forg}(Tj \cdot v).
$$

\end{proof}

\section[Abstract mechanical connection on almost K\"{a}hler manifolds]
{Abstract mechanical connection}

Let $P$ be an almost K\"{a}hler manifold with a complex structure
$\ {\mathcal{J}} : T_x P \rightarrow T_x P \ $, s.t. ${\mathcal{J}}^2 = -1$,
a symplectic form $\omega$ and a ${\mathcal{J}}$-invariant
Riemannian metric ${\mathfrak{s}}$
with the standard relation between these structures
(\ref{akahler})
$$
  \omega (v,w) = 
  {\mathfrak{s}}({\mathcal{J}}v,w) 
  \forall \  v, w \in T_x P.
$$
Let a Lie group $G$ act on $P$ freely and properly 
\footnote{See \cite{Otto87} for some interesting results on how to relax 
the regularity conditions.}
by isometries of the almost K\"{a}hler structure, {\it i.e.} it preserves 
Riemannian, symplectic, and almost complex forms.
The quotient manifold then has a unique Poisson structure such that
the canonical projection $\pi \ : P \rightarrow P/G$ is a Poisson map.
Assume that the $G$ action admits an equivariant momentum map $J$
and that $P/G \stackrel{\pi}\longleftarrow P 
\stackrel{J}\longrightarrow U \subset {\mathfrak{g}}^\ast$
is a full dual pair, {\it i.e.,\/} $\pi$ and $J$ are surjective submersions.
Though we are not interested here in the results for K\"{a}hler
reduction,
\footnote{The reader is refered to, e.g., \cite{HeLo94,Otto87}
for Marsden-Weinstein reduction on  K\"{a}hler manifolds},
we notice that the almost complex structure can be dropped to
the quotient space $P/G$. 
We'll keep the same notation
for the reduced object but we'll write ${\mathcal{J}}(y)$, where
$y=\pi(x)$, to indicate that it can be computed at any 
$x \in \pi^{-1} (y)$.

\subsection{Main constructions}

\begin{dfn}
The \emph{abstract locked inertia tensor} ${\mathbb{I}}(x) :
{\mathfrak{g}} \rightarrow {\mathfrak{g}}^\ast \quad  
\forall \ x \in P$ is defined by the following expression
\begin{equation}
  \label{inert_def}
  \langle {\mathbb{I}} (x) \cdot \xi, \eta \rangle
  = {\mathfrak{s}} (\xi_P(x), \eta_P(x) )
\end{equation}
for any Lie algebra elements $\xi, \eta \in {\mathfrak{g}}$,
where $\xi_P, \eta_P$ are the corresponding infinitesimal generators,
{\it i.e.,\/} vector fields on $P$.
\end{dfn}

The abstract locked inertia tensor is, obviously, an isomorphism for any 
$x \in P$ for which the group action is free. For a general Lie group, 
it is $G$-equivariant in the sense of a map ${\mathbb{I}} \ : P \rightarrow 
\operatorname{L}({\mathfrak{g}}, {\mathfrak{g}}^\ast)$, namely
\begin{equation}
  \label{equiv_I}
{\mathbb{I}}(g \cdot  x) \cdot \operatorname{Ad}_g \xi =
\operatorname{Ad}^\ast_{g^{-1}} {\mathbb{I}}(x) \cdot \xi.
\end{equation}

For an Abelian group, the abstract locked inertia tensor is, in fact,
$G$-invariant and, hence, can be dropped to the quotient $P/G$.
We'll use the same notation
for the reduced object but we'll write ${\mathbb{I}}(y)$, where
$y=\pi(x)$, to indicate that it can be computed at any 
$x \in \pi^{-1} (y)$.

\begin{dfn}
For any choice of a principle connection on $P/G$ define the
\emph{induced metric} ${\mathfrak{s}}'$ on $P/G$ in the following way.
Let $a, b \in T_y (P/G)$ and let $\tilde{a}, \tilde{b}$ be their
corresponding pre-images in the horizontal subspace, {\it i.e.}
$\tilde{a}, \tilde{b} \in \operatorname{hor}(x), \ \pi(\tilde{a}) = a, 
\ \pi(\tilde{a}) = b$, 
where $x \in \pi^{-1}(y)$. As the metric is $G$-invariant we can define
${\mathfrak{s}}' (a,b) ={\mathfrak{s}}_x (\tilde{a},\tilde{b})$,
for any $x \in \pi^{-1}(y)$.
\end{dfn}

\begin{dfn}
The \emph{abstract mechanical connection} on the principle $G$-bundle
$P \rightarrow P/G$ is defined by specifying a horizontal space within 
$T_x P$ at each point 
$x \in P$ to be metric-orthogonal to the tangent to the group orbits:
\begin{equation}
  \label{hor_def}
  \operatorname{hor}(x) = \{ v \in T_x P \ | \ 
  {\mathfrak{s}} (v, \xi_P(x)) = 0
  \quad \forall \ \xi \in {\mathfrak{so}}(3) \}.
\end{equation}
\end{dfn}
The connection one-form ${\mathcal{A}}$  is determined
by $\operatorname{Ker} {\mathcal{A}} (x) 
= \operatorname{hor}(x)$; an  explicit expression for it is
given by the following theorem.

\begin{thm}
Abstract mechanical connection on an almost K\"{a}hler manifold is given by
\begin{equation}
  \label{conn_thm}
  {\mathcal{A}}(x) \cdot w = {\mathbb{I}}^{-1}(x) \cdot 
  {\mathfrak{s}}(\omega^\# (d J(x)), w) \quad \forall \ w \in T_x P.
\end{equation}
\end{thm}

\begin{proof}
For any tangent vector $w \in T_x P$ and any Lie algebra element
$\eta \in {\mathfrak{g}}$ :
$$
{\mathfrak{s}}(w, \eta_P) = {\mathfrak{s}}(w_v + w_h, \eta_P) =
{\mathfrak{s}}(w_v, \eta_P) = {\mathfrak{s}}(\xi^w_P, \eta_P) =
\langle {\mathbb{I}} (x) \xi^w, \eta \rangle,
$$
where $w_v = \xi^w_P$ for some $\xi^w \in {\mathfrak{g}}$ is a vertical
(fiber) component, $w_h$ is a horizontal component, and
${\mathfrak{s}}(w_h, \eta_P) = 0$ by definition.

By definition of the momentum map 
$\eta_P = \omega^\# (d \langle J(x), \eta \rangle )$, so that
$$
{\mathfrak{s}}(w, \eta_P) =  
{\mathfrak{s}}(w, \omega^\# (d \langle J(x), \eta \rangle )) =
\langle {\mathfrak{s}} (\omega^\# (d J), w), \eta \rangle,
$$
where $d J$ is thought of as a ${\mathfrak{g}}^*$-valued one-form on $P$
and we have used the fact that the pairing between ${\mathfrak{g}}$
and its dual is independent of $x \in P$.
Thus,
$$
\langle {\mathbb{I}} (x) \xi^w, \eta \rangle =
\langle {\mathfrak{s}} (\omega^\# (d J), w), \eta \rangle,
$$
and the result follows from the non-degeneracy of the pairing.

To verify that ${\mathcal{A}}$ indeed defines a connection we check
that it satisfies 
${\mathcal{A}}(\xi_P(x)) = \xi, \ \ \forall \ \xi \in {\mathfrak{g}}$ 
and is $G$-equivariant. Consider the pairing of 
${\mathbb{I}}(x) \cdot {\mathcal{A}}(\xi_P(x))$
with an arbitrary element from the Lie algebra $\eta \in {\mathfrak{g}}$
and use above definitions of the connection and 
the abstract locked inertia tensor:
$$
\langle {\mathbb{I}} (x) \cdot {\mathcal{A}}(\xi_P(x)), \eta \rangle =
\langle {\mathfrak{s}}(\omega^\# (d J(x)),\xi_P(x)), \eta \rangle =
{\mathfrak{s}}(\xi_P, \eta_P) = \langle {\mathbb{I}} (x) \xi, \eta \rangle.
$$
From the non-degeneracy of the pairing it follows that
${\mathcal{A}}(\xi_P(x)) = \xi$. The $G$-equivariance means 
$\Phi^\ast_g {\mathcal{A}} = \operatorname{Ad}_g {\mathcal{A}}$
and follows  from equivariance of the momentum map and
equivariance of the abstract locked inertia tensor in the sense
of a map ${\mathbb{I}} \ : P \rightarrow 
\operatorname{L}({\mathfrak{g}}, {\mathfrak{g}}^\ast)$ (see (\ref{equiv_I})).
\end{proof}

\begin{cor} 
\label{cor_A}
The connection one-form can be written as follow :
\begin{equation}
  \label{A_cor}
    {\mathcal{A}}(x) \cdot w = {\mathbb{I}}^{-1}(x) \cdot \operatorname{forg}
    (T J (x) ( {\mathcal{J}} w)) \quad \forall \ w \in T_x P.
\end{equation}
Then,
\begin{equation}
  \label{hor_cor}
  \operatorname{hor}(x) = \operatorname{Ker} (T J(x) \circ {\mathcal{J}}).
\end{equation}
\end{cor}

\begin{proof}
Using (\ref{akahler}), ${\mathcal{J}}^2 = -1$ and omitting $x$
for simplicity we have $\forall \ w \in T_x P$
$$
{\mathcal{A}} \cdot w = {\mathbb{I}}^{-1} \cdot 
{\mathfrak{s}}(\omega^\# (d J), -{\mathcal{J}}^2 w) =
{\mathbb{I}}^{-1} \cdot \omega (\omega^\# (d J), {\mathcal{J}} w) =
{\mathbb{I}}^{-1} \cdot \operatorname{forg} (T J ({\mathcal{J}} w)),
$$
where for the last equality we used the definition of a symplectic form 
and considered the one-form $d_x J$  as a tangent map 
$\operatorname{forg} \circ TJ$ acting on vectors in $T_xP$. 
\end{proof}

\begin{lem}
\label{lem_TJ}
For the choice of the abstract mechanical connection ${\mathcal A}$ on $P$
with $\operatorname{hor} = (\operatorname{Ker} T \pi)^\perp$
the following holds
$$
\operatorname{hor}^\omega = (\operatorname{Ker} T J)^\perp
= {\mathcal J} (\operatorname{Ker} T \pi).
$$
\end{lem}

\begin{proof}
The proof follows readily from (\ref{hor_cor}) of the previous corollary
and the $\omega$-orthogonality of $\operatorname{Ker} T J$ and
$\operatorname{Ker} T \pi$:
\begin{multline}
\nonumber
\operatorname{hor}^\omega 
= (\operatorname{Ker} (T J \circ {\mathcal{J}}))^\omega
= ({\mathcal{J}} (\operatorname{Ker} T J ))^\omega
= \left( ((\operatorname{Ker} T J)^\perp )^\omega \right)^\omega  \\
= (\operatorname{Ker} T J)^\perp 
= ((\operatorname{Ker} T \pi)^\omega)^\perp 
={\mathcal{J}} (\operatorname{Ker} (T J) ),
\end{multline}
where we used that $((W)^\omega)^\perp = 
{\mathcal J} (W)$ for a subspace $W \in T_x P$.
\end{proof}

Below we present two alternative proofs of this lemma
which provide an interesting insight into the issue; these
proofs can be skipped on the first reading.

\emph{Alternative proof.}

By definition, $w \in \operatorname{hor}^\omega (x)$ iff
$\omega(v, w) = {\mathfrak{s}} (v, {\mathcal{J}} w) = 0, \ 
\forall v \in \operatorname{hor} (x)$. Thus, 
$w \in \operatorname{hor}^\omega (x) \Leftrightarrow
{\mathcal{J}} w \in (\operatorname{hor} (x))^\perp$, or 
$$
w \in (\operatorname{hor}^\omega (x))^\perp  \ \Leftrightarrow \
{\mathcal{J}} w \in \operatorname{hor} (x) = 
(\operatorname{Ker} T \pi(x))^\perp.
$$

On the other hand, 
$u \in (\operatorname{Ker} T \pi(x))^\omega =\operatorname{Ker}T J(x)$
iff $\omega (u, \xi_P(x) ) = {\mathfrak{s}} (\xi_P(x),{\mathcal{J}}u) = 0 \ 
\forall \xi \in {\mathfrak{g}}$. Thus, 
$$
u \in \operatorname{Ker}T J(x) 
\Leftrightarrow {\mathcal{J}}u \in (T_x (G \cdot x))^\perp \equiv
(\operatorname{Ker} T \pi(x))^\perp.
$$
Comparing conditions for $u$ and $w$ we conclude that
$\operatorname{hor}^\omega (x) = (\operatorname{Ker}T J(x))^\perp$.
The last equality then follows from 
$\operatorname{Ker}T J = (\operatorname{Ker}T \pi)^\omega$ and
$((\operatorname{Ker}T \pi)^\omega)^\perp 
= {\mathcal J}(\operatorname{Ker}T \pi)$.

\emph{Alternative proof.}

First notice that
$$
\operatorname{hor}^\omega = ((\operatorname{Ker} T \pi)^\perp)^\omega
= ((\operatorname{Ker} T \pi)^\omega)^\perp = 
(\operatorname{Ker} T J)^\perp.
$$
The last equality follows from the following argument. Let $A \subset T_x P$,
then $a \in A^\perp \ \Leftrightarrow \ {\mathfrak s} (a,b) = 0 
\ \forall \ b \in A$. Similarly,
$c \in (A^\perp)^\omega \ \Leftrightarrow \ \omega(c, a) = 0
\ \forall \ a \in A^\perp$. But 
$0 = \omega (c,a) = {\mathfrak s} ({\mathcal J} c, a) 
\ \forall \ a \in A^\perp$ implies that
$$
c \in (A^\perp)^\omega \quad \Leftrightarrow \quad
{\mathcal J} c \in (A^\perp)^\perp \equiv A.
$$
This is equivalent to $c \in {\mathcal J}(A)
\ \Leftrightarrow \ c \in (A^\perp)^\omega$, so that 
$(A^\perp)^\omega = {\mathcal J}(A)$ and
$((\operatorname{Ker} T \pi)^\perp)^\omega = 
{\mathcal J} (\operatorname{Ker} T \pi)$.

Define for any point $\nu \in {\mathfrak{g}}^\ast$   a one-form
${\mathcal{A}}_\nu (x) = \langle \nu, {\mathcal{A}} (x) \rangle$ on $P$.
\begin{lem}
Identifying vectors and one-forms on $P$ via Riemannian metric
$$
({\mathcal{A}}_\nu (x))^\# = ({\mathbb{I}}^{-1}(x) \cdot \nu)_P
$$
\end{lem}
\begin{proof}
Using (\ref{conn_thm}) we obtain
$\forall \ w \in T_x P$:
$$
{\mathcal{A}}_\nu \cdot w = \langle \nu, 
{\mathbb{I}}^{-1} \cdot {\mathfrak{s}} (\omega^\# (d J), w) \rangle =
{\mathfrak{s}} (\omega^\# (d \langle J, {\mathbb{I}}^{-1} \nu \rangle ), w) =
{\mathfrak{s}} (({\mathbb{I}}^{-1} \nu)_P, w).
$$
\end{proof}

\subsection{Abelian groups and reconstruction phases}

In the rest of this section we assume that the Lie group $G$ is Abelian.
A simple corollary of Theorem \ref{thm_bundle} implies that metric orthogonal
horizontal spaces on the Poisson bundle $P \rightarrow P/G$ induce metric 
orthogonal horizontal spaces on symplectic bundles 
$J^{-1}(\mu) \rightarrow P_\mu$ for regular $\mu$.
Hence, by analogy, the reconstruction connection $A$
corresponding to ${\mathcal A}$ by means of Theorem \ref{thm_bundle}
can be called an \emph{abstract mechanical reconstruction connection}.
The same theorem gives us also the corresponding connection on the
symplectic stratification $D = \hat{A}$ by specifying its horizontal
spaces to be $T \pi (\operatorname{hor}^\omega)$.
The following results significantly simplify explicit computations
of these spaces, {\it i.e.} the distribution $D$.

\begin{thm}
For the choice of the abstract mechanical  connection ${\mathcal A}$ 
on $P$, the distribution $D$, which corresponds to the connection on
the symplectic stratification $j : P/G \rightarrow {\mathfrak g}^*$,
is metric orthogonal to the characteristic distribution $E$
in the metric ${\mathfrak{s}}'$ induced on the quotient $P/G$.
Moreover, the distribution $D$ can be explicitly constructed
using the infinitesimal generator vector fields $\xi_P$ according to the
following expression
\begin{equation}
\label{expr_D}
D(y) = T \pi ( {\mathcal J} (\operatorname{Ker} T \pi (x))),
\end{equation}
where $x \in \pi^{-1} (y)$ and
$\operatorname{Ker} T \pi (x) = \operatorname{span}\{ \xi_P (x) \}$.
\end{thm}

\begin{proof}
Consider any vectors $v \in D(y)$ and $w \in E(y) \equiv 
T_y \Sigma_y$. By definition of the induced metric
$
{\mathfrak{s}}' (v, w) = {\mathfrak{s}} (\tilde{v},\tilde{w}),
$
where $\tilde{v}, \tilde{w} \in \operatorname{hor} (x)$ 
are horizontal components
of the pre-images: $T \pi (\tilde{v}) = v, \ T \pi (\tilde{w}) = w$, and
$\ x \in \pi^{-1} (y)$.

From $T \pi (\tilde{v}) = v \in T \pi (\operatorname{hor}^\omega)$ 
it follows that $\tilde{v} \in \operatorname{hor}^\omega 
+ \operatorname{Ker} T \pi$. But $\tilde{v} \in \operatorname{hor} 
\equiv (\operatorname{Ker} T \pi)^\perp$,
so that by Lemma \ref{lem_TJ}
$$
\tilde{v} \in \operatorname{hor}^\omega \cap (\operatorname{Ker} T \pi)^\perp
\equiv (\operatorname{Ker} T J)^\perp \cap (\operatorname{Ker} T \pi)^\perp
$$

For the vector $w \in E(y)$ it holds $T j (w) = 0$ and,
hence, by the commutativity of the diagram in Fig \ref{diag_j}, 
$T J (\tilde{w}) = 0$
for any of its pre-images. In particular,  for the horizontal
pre-image $\tilde{w} \in \operatorname{hor}$ we have
$$
\tilde{w} \in \operatorname{Ker} T J  \cap (\operatorname{Ker} T \pi)^\perp.
$$
From the expressions for $\tilde{v}$ and $\tilde{w}$ it follows that
${\mathfrak{s}}' (v,w) = {\mathfrak{s}} (\tilde{v},\tilde{w}) = 0$.

Finally, (\ref{expr_D}) follows from 
$D = T \pi ((\operatorname{hor})^\omega)$ and Lemma \ref{lem_TJ}.
\end{proof}

\paragraph{Transverse derivatives.}
Here we shall give a new construction of the map $L$ defined by the
Proposition \ref{lem_L} which is crucial for the definition of the 
transverse derivatives, and hence for the computation of the phases .
Our construction is based on Lemma \ref{lem_L_Tj} and
depends implicitly on the choice of the abstract mechanical connection.

\begin{dfn}
For each point $y \in P/G$ define a map 
$N(D, y) : {\mathfrak g} \rightarrow D(y)$ by
\begin{equation}
\label{def_N}
\xi \mapsto T \pi ({\mathcal J}( \xi_P(x))),
\end{equation}
where $\xi$ is a Lie algebra element, $\xi_P(x)$ - its corresponding
infinitesimal generator at $x \in \pi^{-1}(y) \subset P$, and
${\mathcal J}$ is the almost complex structure on $P$.
\end{dfn}

From (\ref{def_N}) it follows that $N$ is a linear map as 
all maps used in its definition are linear.
From the Symplectic Leaf Correspondence Theorem and the fact that $G$
is Abelian and finite dimensional it follows that the dimension of $D(y)$ 
(which equals the co-dimension of the leaf $\Sigma_y$) 
equals the dimension of the algebra ${\mathfrak g}$.
On the other hand, 
$$
\operatorname{dim} (\operatorname{hor}^\omega (x)) =
\operatorname{dim} ({\mathcal J} (\operatorname{Ker} T \pi (x))) =
\operatorname{dim} (\operatorname{Ker} T \pi (x)) =
\operatorname{dim} {\mathfrak g}
$$
Hence, from the fact that $D = T \pi (\operatorname{hor}^\omega)$
the following lemma follows:

\begin{lem}
For each $y \in P/G$ 
the map $N$ is an isomorphism between the Lie algebra ${\mathfrak g}$
and the transverse space $D(y)$ defined at $y$ by
the distribution $D$ on the symplectic
stratification $j : P/G \rightarrow {\mathfrak g}^*$.
\end{lem}


\begin{lem}
\label{lem_com_L}
For an Abelian group $G$, the map $L(D, y)$ defined in the Proposition 
\ref{lem_L} is given by the following composition
\begin{equation}
\label{com_L}
L(D, y) = N(D, y) \circ {\mathbb I}^{-1}(y)
: {\mathfrak g}^* \rightarrow D(y),
\end{equation}
where ${\mathbb I}$ is the abstract locked
inertia tensor.
\end{lem}

\begin{proof}
By the definition of the momentum map, $J_\xi \equiv \langle J, \xi \rangle$
is a Hamiltonian for the vector field $\xi_P$ of the infinitesimal 
transformations, that is, for any vector $u \in T_x P$
$$
\omega(x) (\xi_P, u) = d_x J_\xi (u).
$$
The one-form $d_x J_\xi$ can be thought of as the tangent map $TJ$
acting on vectors in $T_xP$ and paired with $\xi \in {\mathfrak g}$.
Take $u$ to be ${\mathcal J}(\eta_P)$ for some infinitesimal generator
$\eta_P$ corresponding to $\eta \in {\mathfrak g}$. Then,
\begin{multline}
\nonumber
\omega(x) (\xi_P, {\mathcal J}(\eta_P)) = d_x J_\xi ({\mathcal J}(\eta_P)) 
= \langle d_x J ({\mathcal J}(\eta_P)), \xi \rangle 
= \langle \operatorname{forg} (TJ ({\mathcal J}(\eta_P))), \xi \rangle \\
= \langle \operatorname{forg} (Tj \circ T \pi ({\mathcal J}(\eta_P))),
\xi \rangle
= \langle \operatorname{forg} (Tj \circ N (\eta)), \xi \rangle,
\end{multline}
where we used the definition of the map $N$ given by (\ref{def_N}).

On the other hand, 
\begin{multline}
\nonumber
\omega(x) (\xi_P, {\mathcal J}(\eta_P)) =
- \omega(x) ({\mathcal J}(\eta_P),\xi_P) =
{\mathfrak s}(x) (\eta_P,\xi_P) \\
= \langle {\mathbb I}(x) \eta, \xi \rangle
= \langle {\mathbb I}([x]) \eta, \xi \rangle
= \langle {\mathbb I}(y) \eta, \xi \rangle.
\end{multline}

Alternatively, this expression can be obtained from Corollary \ref{cor_A}
using an explicit form of the connection one-form given by (\ref{A_cor}).

From the last two expressions and the non-degeneracy of the pairing
we conclude that $\operatorname{forg} (T j \circ N) = {\mathbb I}$, 
then from Lemma \ref{lem_L_Tj} it follows that
$$
L(D,y) = \operatorname{forg} (T j |_D )^{-1} = 
N(D, y) \circ {\mathbb I}^{-1}(y).
$$
Notice that $L$ is an isomorphism as both $N$ and ${\mathbb I}$ are.
\end{proof}

\paragraph{Dynamic phase.}
Recall that according to (\ref{dyn_ph}), the infinitesimal dynamic
phase is given by the transverse derivative of the reduced Hamiltonian,
which we simplify using the above explicit expression for the map $L$.

\begin{thm}
\label{thm_dp_smpl}
The $\nu$-component of the infinitesimal dynamic phase 
$\xi_{\operatorname{dyn}}$, for any $\nu \in {\mathfrak g}^*$,
can be expressed via the abstract locked inertia tensor and
the almost complex structure according to:
\begin{equation}
\label{dp_smpl}
\langle \nu, \xi_{\operatorname{dyn}}(y) \rangle =
\langle \nu, D_\mu h (y) \rangle =
d h \left( T \pi ({\mathcal J}([x] ) ( ({\mathbb I}^{-1}([x]) 
\cdot \nu)_P )) \right),
\end{equation}
where $x \in [x] = \pi^{-1}(y)$ and  $\mu = j (y)$ 
\end{thm}

\begin{proof}
The proof is quite straightforward and relies on the constructions discussed 
in this section. Using the definition of the transverse derivative,
Lemma \ref{lem_com_L} and $G$-invariance of the abstract
locked inertia tensor and the almost complex structure we obtain
\begin{multline}
\nonumber
\langle \nu, \xi_{\operatorname{dyn}}(y) \rangle =
\langle \nu, D_\mu h (y) \rangle = d h ( L(D,\mu, y) \cdot \nu ) \\
= d h ( N(D, y) \circ {\mathbb I}^{-1}(x) \cdot \nu) = 
d h \left( T \pi ({\mathcal J}(x) ( ({\mathbb I}^{-1}(x) 
\cdot \nu)_P )) \right),
\end{multline}
where the last equality follows from (\ref{def_N}).

As it was pointed out earlier, both ${\mathcal J}$ and 
${\mathbb I}$ are $G$-invariant and, hence, can be dropped
to the quotient $P/G$, so that  (\ref{dp_smpl})  can be
computed at any $x \in [x] \equiv y$. 
\end{proof}

\paragraph{Remark.}
(\ref{dp_smpl}) is equivalent to 
\begin{equation}
\label{dp_smpl2}
\langle \nu, \xi_{\operatorname{dyn}}(y) \rangle =
d H \left( {\mathcal J}(x) ( ({\mathbb I}^{-1}(x)\cdot \nu)_P ) \right),
\end{equation}
where $x \in \pi^{-1}(y) $ and $\mu = j (y)$.

Notice that the (\ref{dp_smpl}) does not depend on the choice of 
$x \in \pi^{-1}(y)$. This agrees with the general 
philosophy of \cite{Blaom99} that all information about the phases
is contained in the reduced quantities. Yet, for the explicit
computations it might be convenient to work with the objects
in the unreduced space. Alternatively, when one has a good model
of the reduced space $P/G$, one can compute a basis ${\mathbf v}_k$ of 
the distribution $D$ at any $y \in P/G$ using isomorphism $L$ and
Lemma \ref{lem_com_L} corresponding to a basis $e_i$ of ${\mathfrak g}^*$.
Then, for any $\nu = \sum \nu^i e_i \in {\mathfrak g}^*$,
the corresponding $\nu$-component of the dynamic phase will be given
by the derivative of the reduced Hamiltonian in the direction
$v = \sum \nu^k {\mathbf v}_k$, {\it i.e.}
$\langle \nu, \xi_{\operatorname{dyn}}(y) \rangle =
d h(\sum \nu^k {\mathbf v}_k)$.

\paragraph{Geometric phase.}
Assuming the $\mu$-regularity of the distribution $D$
(see Section $2.1$), the geometric phase is given by the
transverse derivative of the assembled reduced symplectic form 
$\omega_D$ according to (\ref{omegaD}). In this section we shall
give an explicit construction of this form $\omega_D$ using
the horizontal lifts with respect to the abstract mechanical connection
$A$ and the \emph{unreduced} symplectic form $\omega$. 
This allows us to circumvent explicit computations of the
curvature of the connection one-form that is used in (\ref{geom_ph_curv})
and, in some cases, also the computations of  the
reduced symplectic form that is used in (\ref{omegaD}).

\begin{dfn}
For an Abelian group $G$, define a closed 'horizontal' 
two-form $\omega'$ on $P/G$ according to
\begin{equation}
\label{omega'_def}
\omega' (y) (u, v) := \omega (x) (\tilde{u}, \tilde{v}),
\quad \forall \ u, v \in T_y (P/G), \ y \in P_\mu,
\end{equation}
where $\tilde{u}, \tilde{v} \in {\mathcal A}(x) \equiv \operatorname{hor} (x)
\equiv  (\operatorname{Ker} T \pi )^\perp$ with 
$T \pi (\tilde{u}) = u, T \pi (\tilde{v}) = v$,
$x \in \pi^{-1}(y)$.
\end{dfn}

From the $G$-invariance of the symplectic form $\omega$ as well as
of the horizontal distribution ${\mathcal A}$ we conclude that $\omega'$ 
is well-defined.

\begin{thm}
The two-form $\omega'$ coincides with the assembled two-form $\omega_D$
on $P/G$ : 
$$
\omega' = \omega_D.
$$
\end{thm}

\begin{proof}
We start with the definition of the two-form $\omega'$ above and 
shall demonstrate that the following three special cases hold
for any $y \in P/G$ :
\begin{enumerate}
\item[1)] $\omega' |_E = \omega_\mu \ $, here $E$ is the characteristic
distribution and $\mu = j (y)$,
\item[2)] $\omega' |_D = 0$,
\item[3)] $\omega' (u,v) = 0$ for any $u \in E(y) \equiv T_y P_\mu$
and $v \in D(y)$,
\end{enumerate}
which all together prove the statement of the theorem,
according to the definition of the assembled form (\ref{omegaD}).

\begin{enumerate}
\item[(1)] From the definition of the reduced symplectic form
in the Marsden-Weinstein reduction it follows that
$$
\omega_\mu (y) (u, v) = \omega (x) (\breve{u}, \breve{v}),
$$
where $x \in \pi^{-1}(y) \cap J^{-1}(\mu)$ and
$\breve{u}, \breve{v} \in A(x)$, {\it i.e.} the pre-images
lie in the horizontal space of the reconstruction connection $A$.
Recall that in our case, $A$ denotes the metric orthogonal to the
group orbit within the kernel of $TJ$:
$$
A(x) = (T_x ( G \cdot x))^\perp \cap \operatorname{Ker} TJ (x)
= (\operatorname{Ker} T \pi (x))^\perp \cap \operatorname{Ker} TJ (x).
$$

From Lemma \ref{lem_B1} and the fact that $G$-orbits of any point 
$x \in P$ are isotropic for an Abelian group it follows that
$$
\pi^* E = \operatorname{Ker} T \pi + (\operatorname{Ker} T \pi)^\omega
= (\operatorname{Ker} T \pi)^\omega.
$$

Hence, using the definition on the two-form $\omega'$ above,
for any vectors $u, v \in E = T_y P_\mu$ lying in the
characteristic distribution at $y \in P_\mu$, their pre-images
$\tilde{u}, \tilde{v} \in {\mathcal A}(x)$ satisfy
$$
\tilde{u}, \tilde{v} \in {\mathcal A}(x) \cap 
(\operatorname{Ker} T \pi)^\omega = 
(\operatorname{Ker} T \pi)^\perp \cap \operatorname{Ker} T J
= A(x),
$$
so that
$$
\omega' (y) (u, v) := \omega (x) (\tilde{u}, \tilde{v})
= \omega (x) (\breve{u}, \breve{v}) = \omega_\mu (y) (u, v).
$$

\item[(2)]
Let $u, v \in D(y)$, then $\tilde{u}, \tilde{v} \in \pi^* D = A^\omega$,
by the definition of the reconstruction connection. But 
$\tilde{u}, \tilde{v} \in (\operatorname{Ker} T \pi)^\perp$, so that
\begin{multline}
\nonumber
\tilde{u}, \tilde{v} \in A^\omega \cap (\operatorname{Ker} T \pi)^\perp
= \left( (\operatorname{Ker} T \pi)^\perp \cap \operatorname{Ker} TJ 
\right)^\omega \cap (\operatorname{Ker} T \pi)^\perp \\
= \left( ((\operatorname{Ker} T \pi)^\omega)^\perp 
+ \operatorname{Ker} T \pi \right)
\cap (\operatorname{Ker} T \pi)^\perp.
\end{multline}
Using the modularity property and 
the fact that $\operatorname{Ker} T \pi$ is isotropic, {\it i.e.}
$\operatorname{Ker} T \pi \subset (\operatorname{Ker} T \pi)^\omega$,
and, hence, 
\begin{equation}
\label{isotr}
(\operatorname{Ker} T \pi)^\perp \supset 
((\operatorname{Ker} T \pi)^\omega)^\perp
\end{equation}
we obtain that
$$
\tilde{u}, \tilde{v} \in A^\omega \cap (\operatorname{Ker} T \pi)^\perp
= ((\operatorname{Ker} T \pi)^\omega)^\perp \cap 
(\operatorname{Ker} T \pi)^\perp = ((\operatorname{Ker} T \pi)^\omega)^\perp,
$$
but this space is also isotropic, {\it i.e.} it is contained in its
symplectic orthogonal because of (\ref{isotr}):
$$
((\operatorname{Ker} T \pi)^\omega)^\perp 
= ((\operatorname{Ker} T \pi)^\perp)^\omega \subset 
(((\operatorname{Ker} T \pi)^\omega)^\perp)^\omega.
$$
Thus, $\omega' |_D = 0$.

\item[(3)] 
Finally, combining previous arguments, 
for any $u \in T_y E$ and $v \in D(y)$, $\tilde{u} \in A$ 
and $\tilde{v} \in ((\operatorname{Ker} T \pi)^\omega)^\perp$.
But,
$$
A^\omega = \left( (\operatorname{Ker} T \pi )^\perp 
\cap \operatorname{Ker} TJ \right)^\omega
= ((\operatorname{Ker} T \pi )^\omega)^\perp 
+ (\operatorname{Ker} TJ)^\omega,
$$
so that $\tilde{v} \in A^\omega$ and $\omega (\tilde{u},\tilde{v}) = 0$.
\end{enumerate}

\end{proof}

\begin{cor}
The infinitesimal geometric reconstruction phase is computed
according to 
\begin{equation}
\label{gp_smpl}
\xi_{\operatorname{geom}}(y) = D_\mu \omega' (y),
\end{equation}
where the transverse derivative $D_\mu \omega'$ is computed using 
Lemma \ref{lem_com_L} and the constructions above used for
the computation of the dynamic phase.
\end{cor}


\section{Application: resonant three-wave interaction }

The three-wave equations describe the resonant quadratic nonlinear 
interaction of three waves and are obtained as amplitude equations 
in an asymptotic reduction of primitive equations in optics, fluid 
dynamics and plasma physics. 
It was first analyzed by Alber, Luther, Marsden and Robbins
in \cite{ALMR99} and later in \cite{LAMR00}. 
Here we will only quote the results relevant 
for the definition of the connection and the computation of phases
and refer the reader to \cite{ALMR99} for the  detailed description.
Some results for the Poisson reduction obtained here 
(such expressions for the Casimirs $C_1$ and $C_2$ as well as 
formulas (\ref{red_Pbrack}) for the reduced Poisson bracket
and  (\ref{red_wform}) for the reduced symplectic structure)
are original and were not presented in \cite{ALMR99}.
We shall use the canonical  Hamiltonian structure and ignore an 
alternative Lie-Poisson description of this system.

\paragraph{The phase space and its K\"{a}hler structure.} The phase
space $P$  of the system is ${\mathbb{C}}^3$ with appropriately weighed 
standard 
K\"{a}hler structure. In particular, a $\gamma_i$-weighed canonical Poisson
bracket on ${\mathbb{C}}^3$ is used. This bracket has the real and
imaginary parts of each complex dynamic variable $q_i$ as conjugate
variables. The corresponding symplectic structure is written as follows:
\begin{equation}
  \label{sympl_wave}
  \omega (z, w) = - \sum_k \frac{1}{s_k \gamma_k} 
  \operatorname{Im} (z_k \bar{w}_k),
\end{equation}
where $z, w \in T_q {\mathbb{C}}^3$. 

Similarly, define a weighted  metric on $P$ 
\begin{equation}
  \label{metric_wave}
  {\mathfrak{s}} (z, w) = \sum_k \frac{1}{s_k \gamma_k} 
  \operatorname{Re} (z_k \bar{w}_k),
\end{equation}
and the standard complex structure $ \ {\mathcal{J}} (z) = i z$.
The K\"{a}hler structure  then contains ${\mathfrak{s}}$ and
$\omega$ as real and imaginary parts, respectively.

\paragraph{The symmetry group and momentum map.} Consider the action
of an Abelian group $T^2$ on ${\mathbb{C}}^3$ given by:
 \begin{equation}
   \label{G_action_wave}
(q_1, q_2, q_3) \mapsto 
(\exp^{- i \xi^1} q_1, \exp^{- i (\xi^1+\xi^2)} q_2, \exp^{- i \xi^2} q_3),
 \end{equation}
where $\xi = (\xi_1, \xi_2)$ is an element of the Lie algebra
$t^2 \equiv {\mathbb{R}}^2$.
The vector fields of the infinitesimal transformations corresponding to 
$\xi_1, \xi_2$  are given by:
 \begin{equation}
   \label{infin_gen_wave}
   \xi^1_P (q) = ( -i \xi_1 q_1,  -i \xi_1 q_2, 0), \
   \xi^2_P (q) = ( 0,  -i \xi_2 q_2, -i \xi_2 q_3) \
   \in \ T_q {\mathbb{C}}^3.
 \end{equation}
Points of the form $(q_1, 0,0), (0,q_2,0), (0,0,q_3)$ have non-trivial
isotropy subgroups, and thus account for singularities in the reduced
space, {\it i.e.} as the action is not free, the reduced space fails
to be a smooth manifold (see, e.g. \cite{MaRa1994}).
Henceforth we shall ignore these points and restrict ourself to the
set of regular points in ${\mathbb{C}}^3$.

The momentum map for this action was computed in \cite{ALMR99} 
and is given by 
\begin{equation}
  \label{mom_map_wave}
  J(q_1, q_2, q_3) = (K_1, K_2) = \left( \frac{1}{2} 
\left(\frac{|q_1|^2}{s_1 \gamma_1} + \frac{|q_2|^2}{s_2 \gamma_2} \right),
\frac{1}{2}
\left(\frac{|q_2|^2}{s_2 \gamma_2} + \frac{|q_3|^2}{s_3 \gamma_3} \right)
\right).
\end{equation}
We shall keep the notations $(K_1, K_2)$ for the values of the momentum map
to be consistent with  \cite{MaRa1994}; they play the role
of the $\mu = J(x)$ in the  exposition above.

It is checked directly using (\ref{G_action_wave}) that the
momentum map is $G$-invariant. 
For further applications we note that even though $J$ is not analytic,
we can consider its differential as a real-valued map of the tangent space
$T {\mathbb{C}}^3$ to $T {\mathbb{R}}^2$.

\paragraph{The Hamiltonian.} The Hamiltonian for the three-wave interaction
is
$$
H = -\frac{1}{2} (\bar{q}_1 q_2 \bar{q}_3 + q_1 \bar{q}_2 q_3 ).
$$
Hamilton's equations are $\dot{q_k} = \{ q_k , H \}$ and it is
straightforward to check that in complex notations they are given by
$$
\frac{d q_k}{d t}=- 2 i s_k \gamma_k \frac{\partial H}{\partial \bar{q}_k}.
$$

\paragraph{Poisson reduction.} 
It was shown in \cite{ALMR99} that the following quantities constitute 
invariants for the $T^2$ action
$$
X+i Y = q_1 \bar{q}_2 q_3, \quad Z_1 = |q_1|^2 - |q_2|^2, 
\quad Z_2 = |q_2|^2 - |q_3|^2.
$$
They provide coordinates for the four-dimensional (real) orbit space.
The symplectic leaves in it are two-dimensional.
This follows from the Leaf Correspondence Theorem, as $T^2$ is 
Abelian, and each point in it being a coadjoint orbit has co-dimension $2$.
One can define two Casimirs on ${\mathbb{C}}^3/T^2$, e.g.
$$
C_1 = (X^2 + Y^2) - \kappa_4 (2 s_2 \gamma_2 K_1 + Z_1) 
(2 s_3 \gamma_3 K_2 + Z_2) (2 s_2 \gamma_2 K_2 - Z_2) ],
$$
where $\kappa_4 = (s_1 \gamma_1 s_2 \gamma_2 s_3 \gamma_3)/
(s_1 \gamma_1 + s_2 \gamma_2)(s_2 \gamma_2 + s_3 \gamma_3)^2$, and
$$
C_2 =  (Z_1 - 2 s_1 \gamma_1 K_1) (s_2 \gamma_2 + s_3 \gamma_3) + 
(Z_2 + 2 s_3 \gamma_3 K_2) (s_1 \gamma_1 + s_2 \gamma_2),
$$
which can be obtained by a pull-back of properly defined 
Casimirs on ${\mathfrak{g}}^\ast$
using ideas of dual pairs (see, e.g. \cite{Weinstein83, Blaom1998}).
The level set of these Casimirs, defined for any momentum map value
$(K_1, K_2) \in U$, where $U \equiv J({\mathbb C}^3) \subset t^2$,
by the set $\{ C_1 = 0, C_2 = 0 \}$ determines the corresponding symplectic 
leaf $P_{(K_1, K_2)}$ in the reduced space.

The reduced Poisson bracket on ${\mathbb C}^3/T^2 \cong (X, Y, Z_1, Z_2)$
is given for any two functions $f, k$ by
\begin{equation}
  \label{red_Pbrack}
\{ f, k \} = \det ( \nabla C_2 \ \nabla C_1  \ \nabla f \ \nabla k \ ).
\end{equation}
The reduced Hamiltonian equations of motion have the following form
$$
\dot{X} = 0, \
\dot{Y} = -\frac{\partial C_2}{\partial Z_2}\frac{\partial C_1}{\partial Z_1}
+ \frac{\partial C_2}{\partial Z_1}\frac{\partial C_1}{\partial Z_2}, \
\dot{Z}_1 = \frac{\partial C_2}{\partial Z_2}\frac{\partial C_1}{\partial Y},\
\dot{Z}_1 = -\frac{\partial C_2}{\partial Z_1}\frac{\partial C_1}{\partial Y}.
$$

Notice that the second Casimir $C_2$ establishes a linear dependence between
$Z_1$ and $Z_2$; hence, we can solve for one of them, say $Z_1$ and restrict
ourselves to the consideration of three-dimensional subspace in 
${\mathbb{C}}^3/T^2$ defined by $(X, Y, Z_2)$. The dynamics in $Z_1$
can then be trivially reconstructed. In this case, the first Casimir $C_1$ 
can be rewritten as
\begin{equation}
  \label{phi_def_wave}
  \phi = (s_2 \gamma_2 + s_3 \gamma_3) [ (X^2 + Y^2) 
  - \kappa_3 (\delta - Z_2) (2 s_3 \gamma_3 K_2 + Z_2)
    (2 s_2 \gamma_2 K_2 - Z_2) ],
\end{equation}
where $\kappa_3 = (s_1 \gamma_1 s_2 \gamma_2 s_3 \gamma_3)/
(s_2 \gamma_2 + s_3 \gamma_3)^3$ and 
$\delta = 2 s_2 \gamma_2 K_1 + 2 s_3 \gamma_3 (K_1 - K_2)$.
This relation defines two-dimensional (perhaps singular) surfaces
in $(X, Y, Z_2)$ space, with $Z_1$ determined by the values of the
invariants and conserved quantities. These surfaces are called
\emph{three-wave surfaces}. 

The reduced Poisson bracket in $(X, Y, Z_2)$ space is given by
$$
\{ f, k \} = \nabla \phi \cdot (\nabla f \times \nabla k),
$$
for any functions $f, k$.
For a non-singular point $y$ on a symplectic leaf $P_{(K_1, K_2)}$
the induced symplectic form is then given by
\begin{equation}
  \label{red_wform}
\omega_{(K_1, K_2)} (v,w) = 
- \frac{\nabla \phi}{\|\nabla \phi \|^2}\cdot (v \times w),
\end{equation}
where $v, w \in T_y P_{(K_1, K_2)}$; here $(K_1, K_2)$ is the momentum
value which determines a particular symplectic leaf $P_{(K_1, K_2)}$. 
Thus, for a function $f$ on the orbit space  $(X, Y, Z_2)$, the 
corresponding Hamiltonian vector field has the form
$X_f = - \nabla \phi \times \nabla f$.

The reduced Hamiltonian for the three-wave interaction is given by
$h = - X$ and produces the following reduced equations of motion
$$
\dot{X} = 0, \quad \dot{Y} = \frac{\partial \phi}{\partial Z_2},
\quad \dot{Z}_2 = -2 (s_2 \gamma_2 + s_3 \gamma_3) Y,
$$
which otherwise can be obtained by the restriction of the equations
of motion in ${\mathbb C}^3/T^2 \cong (X, Y, Z_1, Z_2)$ to three-wave
surfaces.

\paragraph{Abstract Mechanical connection.} 
First of all, we compute the locked inertia tensor 
${\mathbb{I}}(q) \ : \ {\mathbb{R}}^2 \mapsto {\mathbb{R}}^2$
using its definition (\ref{inert_def}).
It is an isomorphism for regular points $q$ and is given by the following
expression
\begin{equation}
  \label{inert_wave}
  {\mathbb{I}}(q) = 
\begin{pmatrix}
  2 K_1 \ & \dfrac{|q_2|^2}{s_2 \gamma_2} \\
        &                              \\
  \dfrac{|q_2|^2}{s_2 \gamma_2} & 2 K_2 \
\end{pmatrix},
\end{equation}
where $K_1(q_1, q_2, q_3)$ and $K_2(q_1, q_2, q_3)$ are components of 
the momentum map given by (\ref{mom_map_wave}).
Notice that ${\mathbb{I}}$ can be dropped to the quotient space:
\begin{equation}
  \label{inert_wave2}
  {\mathbb{I}}(X, Y, Z_1, Z_2) = 
\begin{pmatrix}
  2 K_1 \ & \dfrac{2 s_3 \gamma_3 K_2 + Z_2}{s_2 \gamma_2 + s_3 \gamma_3} \\
        &                              \\
  \dfrac{2 s_3 \gamma_3 K_2 + Z_2}{s_2 \gamma_2 + s_3 \gamma_3} & 2 K_2 \
\end{pmatrix},
\end{equation}
where $K_1, K_2$ are now functions of $(X, Y, Z_1, Z_2)$ as
the momentum map factors through the quotient and are constant
on each symplectic leaf, or a three-wave surface of thereof, in 
the reduced space.

Using Corollary \ref{A_cor} we can explicitly construct the corresponding
abstract mechanical connection one-form 
$$
{\mathcal{A}}(x) \cdot w = 
 - {\mathbb{I}}^{-1}(x) \cdot  \operatorname{Im} 
\begin{pmatrix}
\dfrac{\bar{q}_1 w_1}{s_1 \gamma_1}+\dfrac{\bar{q}_2 w_2}{s_2 \gamma_2} \\
\\
\dfrac{\bar{q}_2 w_2}{s_2 \gamma_2}+\dfrac{\bar{q}_3 w_3}{s_3 \gamma_3}
\end{pmatrix}.
$$
Comparing the last expression with the definition of the 
Riemannian structure (\ref{metric_wave}) we can conclude that 
$$
w \in \operatorname{hor} (q) 
\ \Leftrightarrow \
\left\{ \begin{matrix}
\operatorname{Im} \left(
\dfrac{\bar{q}_1 w_1}{s_1 \gamma_1} + \dfrac{\bar{q}_2 w_2}{s_2 \gamma_2}
\right) = 0  \\
\\
\operatorname{Im} \left(
\dfrac{\bar{q}_2 w_2}{s_2 \gamma_2} + \dfrac{\bar{q}_3 w_3}{s_3 \gamma_3}
\right) = 0 
\end{matrix} \right .
\ \Leftrightarrow \
\left\{ \begin{matrix}
{\mathfrak{s}}(\xi^1_P(q), w) = 0 \\
{\mathfrak{s}}(\xi^2_P(q), w) = 0 
\end{matrix} \right . ,
$$
that is, the horizontal space is precisely the metric orthogonal
to the group orbits.
Similar computations show that $(\operatorname{hor}(q))^\omega$ 
is determined by the span of vectors 
${\mathcal{J}}(\xi^1_P(q)), {\mathcal{J}}(\xi^2_P(q))$, and
$$
(\operatorname{hor}(q))^\omega =
\operatorname{span} \{ i \xi^1_P(q), i \xi^2_P(q) \} = 
(\operatorname{Ker} T J(q))^\perp.
$$

The distribution $D = \hat{A}$ on the symplectic stratification 
$j  : \, P/G \rightarrow {\mathfrak{g}}^\ast$ is obtained
by applying the tangent map $T \pi$ to the space 
$(\operatorname{hor}(q))^\omega$. Define vectors 
${\mathbf v}_1, {\mathbf v}_2$ tangent to the quotient space at the point
$(X, Y, Z_1, Z_2)$ to be the images of a basis in 
$(\operatorname{hor}(q))^\omega$ under this tangent map:
$$
{\mathbf v}_1 = T \pi (q) (i \xi^1_P) = 2 \xi^1 \cdot (X, Y, Z_1, 
\frac{s_2 \gamma_2}{s_2 \gamma_2 + s_3 \gamma_3}(2 s_3 \gamma_3 K_2 + Z_2) )
$$
and
$$
{\mathbf v}_2 = T \pi (q) (i \xi^2_P) = 2 \xi^2 \cdot (X, Y, 
-\frac{s_2 \gamma_2}{s_1 \gamma_1 + s_2 \gamma_2}(2 s_1 \gamma_1 K_1 - Z_1),
Z_2).
$$
Then,
\begin{equation}
  \label{hor_D_wave}
  D (X, Y, Z_1, Z_2) = T \pi (\operatorname{hor}(q))^\omega  =
  \operatorname{span} \{ {\mathbf v}_1, {\mathbf v}_2 \}. 
\end{equation}

To finish the construction of the map $L$ given by Lemma \ref{lem_com_L}
we need to substitute 
${\mathbb I}^{-1} \cdot \nu, \ \forall \nu \in (t^2)^*$ for the $\xi$. 
Then, for any $\nu = \sum \nu_k e^k \in (t^2)^*$, where
$e^k$ is the dual basis of $(t^2)^*$, the map $L$ is given by
$$
L : \nu \mapsto L(\nu) = \sum \nu_k T \pi (i ({\mathbb I}^{-1} \cdot \nu)^k_P)
= \sum \nu_k {\mathbf v}_k,
$$
where in the expressions for ${\mathbf v}_k$ 
we take $\xi^k = ({\mathbb I}^{-1} \cdot \nu)^k$.



\paragraph{Phases for the three-wave interaction.} 
Recall that the reduced Hamiltonian on $P/G$ is $h = - X$.
Applying Theorem \ref{thm_dp_smpl} we can immediately obtain
the $\nu$-component the associated dynamic phase by computing
directional derivatives of the reduced Hamiltonian in the
directions ${\mathbf v}={\mathbf v}_1 + {\mathbf v}_2$ 
in the transverse distribution $D$:
$$
\langle \nu , \xi_{\operatorname{dyn}} \rangle =
d h ({\mathbf v}) 
= \frac{2 h}{\det {\mathbb I}} \left( 2 K_1 \nu_2 + 2 K_2 \nu_1 - 
\frac{2 s_3 \gamma_3 K_2 + Z_2}{s_2 \gamma_2 + s_3 \gamma_3} (\nu_1 + \nu_2)
\right),
$$
where $(K_1, K_2)$ are the momentum values at $(X, Y, Z_1, Z_2)$
along the reduced trajectory $y_t$. To get the dynamic phase
$g_{\operatorname{dyn}}$ one integrates the exponent of this expression
along the reduced trajectory $y_t$ on a three-wave surface.

The infinitesimal geometric phase $\xi_{\operatorname{geom}}$,
as a two-form on the reduced space, can be computed using
(\ref{gp_smpl}), so that  its $\nu$-component is given by
$$
\langle \nu ,\xi_{\operatorname{geom}}(y) \rangle
= \langle \nu , D_\mu \omega' (y)\rangle.
$$
These expression should be computed using standard formulas for
the differentials of $p$-forms. We omit here the calculations
of the dynamic phase 
as they crucially depend on the area over which the two-form
is integrated.

\section{Concluding remarks}

If the phase space $P$ has an almost K\"{a}hler structure
a preferred connection can be defined
by declaring horizontal spaces at each point to be metric orthogonal to 
the tangent to the group orbit. We call it \emph{abstract mechanical
 connection}. 
Then, explicit formulas  for the corresponding ${\mathfrak{g}}^\ast$-valued
one-form ${\mathcal{A}}$ in terms of the momentum map, symplectic and
complex structures can be derived. 
Also, we show that in this case the horizontal spaces for the 
induced connections are metric orthogonal to the corresponding
natural vertical spaces for each foliation.

These results are applied to the resonant three-wave interaction problem
(see, e.g., \cite{ALMR99}). The corresponding horizontal spaces are constructed
and a formula for the dynamic phase  is obtained.
The associated geometric phase is given by the integral
of a two-form which is defined by the reduced symplectic structure.

\section*{Acknowledgments}
The authors would like to thank Tudor Ratiu and Sameer 
Jalnapurkar for helpful comments.



\end{document}